\documentclass[Euromicro PDP 2025]{IEEEtran}
\IEEEoverridecommandlockouts
\usepackage{algorithm}
\usepackage{algpseudocode}
\usepackage{graphicx}
\usepackage{float}
\usepackage{hyperref}
\usepackage{color}
\usepackage{amsmath}
\usepackage{mathtools}
\newcommand{\R}{{\mathcal R}}
    
\begin{document}

\title{Communication-reduced Conjugate Gradient Variants for GPU-accelerated Clusters
\thanks{This work was partially supported by: Spoke 6 ``Multiscale Modelling \& Engineering Applications” of the Italian Research
 Center on High-Performance Computing, Big Data and Quantum Computing (ICSC) funded by
 MUR-NextGenerationEU (NGEU); the ``Energy Oriented Center
 of Excellence (EoCoE III): Fostering the European Energy Transition with Exascale” EuroHPC
 Project N. 101144014, funded by European Commission (EC).}
}

\author{\IEEEauthorblockN{Massimo Bernaschi, Mauro G. Carrozzo, Alessandro Celestini, Giacomo Piperno}\\
\IEEEauthorblockA{\textit{Institute for Applied Computing (IAC), 
Via dei Taurini, 19, I-00185 Rome\\
name.surname@cnr.it}}\\
\and
\IEEEauthorblockN{Pasqua D'Ambra}\\
\IEEEauthorblockA{\textit{Institute for Applied Computing (IAC), 
Via P. Castellino, 111, I-80131 Naples\\
pasqua.dambra@cnr.it}
}}

\maketitle

\begin{abstract}
Linear solvers are key components in any software platform for scientific and engineering computing. The solution of large and sparse linear systems lies at the core of physics-driven numerical simulations relying on partial differential equations (PDEs) and often represents a significant bottleneck in data-driven procedures, such as scientific machine learning. In this paper, we present an efficient implementation of the preconditioned s-step Conjugate Gradient (CG) method, originally proposed by Chronopoulos and Gear in 1989, for large clusters of Nvidia GPU-accelerated computing nodes.

The method, often referred to as communication-reduced or communication-avoiding CG, reduces global synchronizations and data communication steps compared to the standard approach, enhancing strong and weak scalability on parallel computers. Our main contribution is the design of a parallel solver that fully exploits the aggregation of low-granularity operations inherent to the s-step CG method to leverage the high throughput of GPU accelerators. Additionally, it applies overlap between data communication and computation in the multi-GPU sparse matrix-vector product.

Experiments on classic benchmark datasets, derived from the discretization of the Poisson PDE, demonstrate the potential of the method.
\end{abstract}


\section{Introduction}
\label{intro}

Solving large, sparse linear systems of algebraic equations of the form:
\begin{equation}
Ax=b, \  A \in \R^{n \times n} \ \ \  x, b \in \R^{n}
\label{sys}
\end{equation}
is fundamental in numerous computational models and simulation software utilized in Computational and Data Science. Often, the matrix $A$ is symmetric positive definite (spd), as discussed in~\cite{MR3425300} and a common approach to the solution is a preconditioned version of the Conjugate Gradient (CG) method, which belongs to the general class of Krylov Subspace projections Methods (KSMs)~\cite{saadbook}. A KSM is an iterative process that extracts an approximate solution to problem (\ref{sys}) from a Krylov vector subspace of the form:
\[
K_m(v)= span\{ v, Av, \ldots, A^{m-1}v \}, 
\]
where $A$ is the matrix in (\ref{sys}) and $v \in \R^{n}-\{0\}$.
The subspace $K_m$ serves as the search space for potential solutions and has dimension $m$; thus, $m$ constraints, known as orthogonality conditions, are required to determine a unique solution. In the CG method, solution $x$ belongs to the affine vector space $x_0+K_m(r_0)$, where $r_0=b-Ax_0$ and $x_0$ is the starting solution. The orthogonality conditions
read as $r_k=b-Ax_k \perp K_m, \ \forall k=1, 2,  \ldots $. It can be demonstrated that the solution $x$ obtained by the CG method minimizes the A-norm\footnote{If $A$ is spd, the A-inner product is defined as $(x,y)_A=x^TAy$ and the A-norm of $x \in \R^n$ is $\|x\|_A=\sqrt{x^TAx}$.} of the error $E(x)=\|x-\overline{x}\|_A, \ \forall x \in K_m(r_0)$, where $\overline{x}=A^{-1}b$ is the exact solution. We finally observe that vectors in $K_m(r_0)$, i.e., the CG search directions, are chosen to be A-orthogonals or A-conjugates, i.e., they are orthogonal with respect to the A-inner product. This selection contributes to the favorable convergence properties of the CG method when $A$
is well-conditioned, i.e., when the ratio of its largest to smallest eigenvalue is small. For more general spd matrices, with potentially challenging eigenvalue distributions (ill-conditioned matrices), enhancing CG convergence is achieved by employing a preconditioned version of CG. This involves solving a modified system, such as $B^{-1}Ax=B^{-1}b$, where $B$ is an spd preconditioner. The preconditioner 
$B$ is chosen such that its inverse is computationally inexpensive to apply (ideally much less costly than inverting $A$) and 
$B^{-1}A$ has a more favorable condition number than $A$.

Efficient implementation of the preconditioned CG method on modern clusters with heterogeneous computing nodes—featuring multiple CPUs and GPU accelerators that achieve high performance through various forms of parallelism and deep memory hierarchies—remains a significant challenge. This is evidenced by the relatively small percentage of peak performance achieved by supercomputers when measured with the HPCG benchmark~\cite{HPCG}.
The primary challenge arises from the low computation-to-data movement ratio of the basic operations involved in the method, such as sparse matrix-vector products and preconditioner applications, as well as the frequent global synchronizations required for scalar products and norm computations at each iteration. The HPCG benchmark was specifically introduced to evaluate computing systems under these memory and communication-bound computations, which involve irregular data access patterns and fine-grained operations. It is based on the standard formulation of the CG method coupled with a multigrid preconditioner, whose application requires inherently sequential triangular solvers. In response to these limitations, a complementary approach has been pursued by the numerical software community over the years. Significant efforts have been directed toward developing variants of the CG algorithm and preconditioners that reduce global synchronizations and minimize or hide communication overheads. These advancements in algorithm design aim to mitigate the increasing disparity between data movement speeds and computational speeds in modern computing architectures (see Fig.2 in~\cite{Dongarra2022}). 

This paper presents an efficient implementation of a well-known reduced-communication variant of the preconditioned CG method, referred to as the preconditioned s-step CG, for large heterogeneous clusters equipped with the latest-generation Nvidia GPUs. In the literature, the s-step CG method is also commonly known as Communication-Avoiding CG (CA-CG).
To the best of our knowledge, our software represents the first publicly available open-source implementation of this method for heterogeneous clusters. It is part of a broader research effort aimed at developing a general-purpose library of preconditioners and solvers tailored for distributed multi-GPU systems. This paper focuses on the performance of the multi-GPU s-step CG method, combined with a parallel preconditioner whose efficient implementation is included in the library.
The main contributions of this paper are as follows:
i) design and implementation of fundamental computational building blocks for a class of communication-reduced Krylov Subspace Methods (KSMs) that exploit computation aggregation and efficient data access patterns in a multi-GPU environment;
ii) development of a freely available, modular, and extensible software library, named {\em BootCMatchGX}, featuring communication-reduced sparse linear solvers and preconditioners designed for clusters of GPU-accelerated computing nodes;
iii) performance analysis of a hybrid MPI-CUDA implementation of the preconditioned s-step CG method.

The paper is organized as follows. In section \ref{background} we describe the classic CG method and main principles leading to variants which reduce communication and global synchronization. Detailed description of the preconditioned $s$-step CG and its computational building blocks and complexity is in Section \ref{psstepcg}. Section \ref{design} describes design principle and implementation issues of main kernels, while section \ref{prec} describes key features of the parallel preconditioner we apply in this work. Section \ref{relwork} summarizes related works. 
In Section \ref{results} we discuss performance results obtained on a cluster of GPUs for solving systems of size up to $1$ billion.
Some concluding remarks and future work are outlined in Section \ref{conclusions}. 

\section{Background}
\label{background}

In the following we present the classic algorithm originally derived for the CG method introduced by Hestenes and Stiefel~\cite{HS1952} and some of the variants aimed at reducing the number of global synchronization steps needed for scalar (dot) products of vectors involved in orthogonalization procedures, on large-scale parallel machines. It is worth noting that, although the various formulations of the method are equivalent from the mathematical point of view, their behaviour in inexact arithmetic is different, indeed in some cases, numerical instabilities affect maximum attainable accuracy and convergence, as widely discussed in the literature~\cite{carsonetal2018,GLC2021}.   
 Algorithm \ref{alg:pcghs1} shows the preconditioned CG applied to a linear system as in (\ref{sys}), where $B \in \mathcal{R}^{n \times n}$ is a spd preconditioner. 
{\small 
\begin{algorithm}
  \caption{$PCG-HS(A,b,B,x_0,rtol,kmax)$}
  \label{alg:pcghs1}
  \begin{algorithmic}[1]
  \State $r_0=b-A x_0$ 
  \State $u_0=B^{-1} r_0$ 
  \State $p_0=u_0$
    \For{$k=0,1,...,kmax$}
    	\State $s_k=Ap_k$
             \State $\alpha_k=u_k^Tr_k/p_k^Ts_k$
             \State $x_{k+1}=x_k+\alpha_k p_k$
             \State $r_{k+1}=r_k-\alpha_k s_k$
             \State $u_{k+1}=B^{-1}r_{k+1}$
             \State $\beta_k=u_{k+1}^Tr_{k+1}/u_k^Tr_k$
             \State $p_{k+1}=u_{k+1}+\beta_k p_k$
    	\If{($\|r_{k+1}\| < rtol$)}
    	\State return $x_{k+1}$
	    \EndIf
    \EndFor
    \State Return error, algorithm did not converge within $kmax$ iterations
  \end{algorithmic}
\end{algorithm}
}

We observe that $\alpha_k$ and $\beta_k$ computations in the previous algorithm require two synchronization points because $\beta_k$ can not be computed if the residual vector $r_k$ and the preconditioned residual $u_k$ have not been updated. A first attempt to reduce synchronization points was introduced by Saad in~\cite{Saad1985}. Since the two vectors $r_{k+1}$ and $r_k$ are B-orthogonal, i.e. are orthogonal with respect to the B-inner product, we can write:
$\|r_{k+1} \|_B^2 + \|r_k \|_B^2 = \alpha_k^2 \|s_k\|_B^2$,
then, we can obtain $r_{k+1}^Tu_{k+1}$ by $r_k^Tu_k$ and $s_k$, and it can be computed simultaneously to $p_k^Ts_k$. That is we have the following formula for $\beta_k$:
\[
\beta_k= \alpha_k^2 \frac{\|s_k\|_B^2}{\|r_k\|_B^2}-1.
\]
Saad emphasized that this approach is unstable; it should be used together with possible re-computation of $\beta_k$ as soon as $r_{k+1}$ is available at the expense of an additional dot product (see~\cite{Meurant1987}). Similar arguments to reduce global synchronizations have been first used in~\cite{DR1992} for a version of CG where the dot product $p_k^Ts_k$ is avoided by using alternative expressions involving already computed dot products. Although no theoretical proof has been shown, the authors of~\cite{DR1992} did many experiments on problems with critical eigenvalue distributions (ill-conditioned systems) and no instabilities appear. Some years later, Notay applied a similar approach to have a variant of a flexible version\footnote{A flexible version of a KSM is a variant that preserves convergence property when the preconditioner $B$ is an operator varying along the iterations.} of the preconditioned CG (see~\cite{NN2015}) requiring only one synchronization point. In this flexible version, at each iteration, computational building blocks are: one sparse matrix-vector product (SpMV), one preconditioner application (represented by the SpMV involving $B^{-1}$), three coupled dot products, which can be computed by a single global synchronization and exploiting data re-use at the GPU-level to optimize memory accesses, and finally two couple of axpy operations, which can be computed exploiting again data re-use. 
Note that the additional dot product and axpy operation, with respect to Algorithm \ref{alg:pcghs1}, are generally amortized in a parallel setting by the reduction in synchronization steps.
This last algorithm has been implemented in the {\em BootCMatchGX} library~\cite{parco2020,BERNASCHI2020100041,ieee2023}. 
{\small 
\begin{algorithm}
\caption{$PFCG-NN(A,b,B,x_0,rtol,kmax)$}
  \label{alg:pfcgnn}
  \begin{algorithmic}[1]
\State $r_0=b-A x_0$ 
\State $u_0 = p_0 = B^{-1} r_0$ \label{line:prec0}
\State $w_0 = s_0 = A u_0$
\State $\alpha_0 = u_0^T r_0$
\State $\beta_0 = \rho_0 = p_0^T s_0$
\State $x_{1} = x_0 + \alpha_0 / \rho_0 p_0$
\State $r_{1} = r_0 - \alpha_0 / \rho_0 s_0$
\For{$k=1,...,kmax$}
	\State $u_{k} = B^{-1} r_{k}$ 
	\State $w_{k} = A u_{k}$ 
	\State $\alpha_k = u_{k}^T r_{k}$ 
	\State $\beta_k = u_{k}^T w_{k}$
	\State $\gamma_k = u_{k}^T s_{k-1}$
	\State $\rho_k = \beta_k - \gamma_k^2 / \rho_{k-1}$
	\State $p_{k} = u_{k} - \gamma_k / \rho_{k-1} p_{k-1}$
	\State $x_{k+1} = x_k + \alpha_k / \rho_k p_{k}$
	\State $s_{k} = w_{k} - \gamma_k / \rho_{k-1} s_{k-1}$
	\State $r_{k+1} = r_k - \alpha_k / \rho_k s_{k}$ 
        \If{($\|r_{k+1}\| < rtol$)}
    	\State return $x_k$
	\EndIf
\EndFor
\State Return error, algorithm did not converge within $kmax$ iterations
\end{algorithmic}
\end{algorithm}
}

\subsection{Preconditioned $s$-step CG}
\label{psstepcg}

The goal of reducing synchronization points for dot product computations motivated Chronopoulos and Gear to propose the so-called $s$-step CG method. This method groups the vector operations from $s$ consecutive iterations into a single, larger block iteration. By doing so, for a fixed number of iterations, the approach theoretically allows for an $\mathcal{O}(s)$ reduction in communication and synchronization costs. Additionally, this method increases the granularity of basic computations, which can improve the efficiency of high-throughput processors by better exploiting their capabilities.

The implementation of the preconditioned $s$-step method, as described in ~\cite{ChronGear1989p}, utilizes $s$ linearly independent A-conjugate directions. These directions are defined as: $
K_s(r_i)=\{ B^{-1}r_i, \ldots, (B^{-1}A)^{s-1}r_i \}$.
These directions are then used to compute the new $s$-dimensional solution at iteration $i$ within the Krylov subspace: $
    K_{is+1}(r_0)=span \{ B^{-1}r_0, B^{-1}Ar_0, \ldots, (B^{-1}A)^{is} r_0 \}$
The newly generated directions, which can be identified as $p_i^1, \ldots, p_i^s$, must be made A-conjugate to the previous $s$ directions. Finally, the error function $E(x) = |x - \overline{x}|_A^2$ is minimized with respect to all the new $s$ directions to obtain the updated solution $x_{i+1}$ and the new residual $r_{i+1} = b - Ax_{i+1}$.

For completeness, we summarize the method in detail as described in ~\cite{ChronGear1989p,ChronGear1989}. Let the new solution $x_{i+1}$ in the affine vector space $L_i = { x_i + \sum_{j=1}^s \alpha_i^j p_i^j }$ be expressed as: $x_{i+1}=x_i+\alpha_i^1p_i^1+\alpha_i^2p_i^2+\ldots+\alpha_i^sp_i^s$.
The coefficients $\alpha_i^j$ can be obtained by solving the following optimization problem:
\begin{equation}
\label{min}
\alpha_i^j=\arg \min_{x \in L_i} E(x).
\end{equation}
It can be shown that the residual vector $r_{i+1} = b - Ax_{i+1} = r_i - \alpha_i^1 A p_i^1 - \alpha_i^2 A p_i^2 - \ldots - \alpha_i^s A p_i^s$ is the gradient of $E(x)$ at $x_{i+1}$. Therefore, solving problem (\ref{min}) is equivalent to enforcing that $r_{i+1}$ is orthogonal to $L_i$.
Let $W_i = { (p_i^j)^T A p_i^l }, \ 1 \leq j,l \leq s$ be the $s \times s$ non-singular spd matrix formed by the dot products involving the linearly independent vectors $p_i^1, \ldots, p_i^s$. Let $\alpha_i = (\alpha_i^1, \ldots, \alpha_i^s)^T$ represent the $s$-dimensional vector of coefficients. These coefficients can be obtained by solving the system of normal equations:
\begin{equation}
\label{normaleq}
W_i \alpha_i = m_i,
\end{equation}
where $m_i=(r_i^Tp_i^1, \ldots, r_i^Tp_i^s)^T$.
To compute the directions $p_i^1, \dots, p_i^s$ from the previous ones, the A-conjugation property can be utilized. Specifically, we express the directions as follows:
$p_{i}^j=r_{i}+\beta_{i-1}^{(j,1)}p_{i-1}^1+\ldots+\beta_{i-1}^{(j,s)}p_{i-1}^s \ \ j=1, \ldots, s$.
We then impose the following conditions:
\[
W_{i-1} \beta_{i-1}^j+c_{i}^j = 0 \ \ j=1, \ldots, s,
\]
where $\beta_{i-1}^j = (\beta_{i-1}^{(j,1)}, \ldots, \beta_{i-1}^{(j,s)})^T$ and $c_i^j=(r_{i}^TA^jp_{i-1}^1, \ldots, r_{i}^TA^jp_{i-1}^s)^T, \ j=1, \ldots, s$. 
Finally, to compute the coefficients $\beta_{i-1}^j$ and $\alpha_i$, it is necessary to solve $s+1$ spd systems of dimension $s$, which involve the matrices $W_i$.
In ~\cite{ChronGear1989p, ChronGear1989}, recurrence formulas are derived to compute the coefficients $c_i^j$ (see Proposition 4.2 in ~\cite{ChronGear1989}), the matrix $W_i$ (see Corollary 4.2 in ~\cite{ChronGear1989}), and the vectors $m_i$ (see Corollary 4.3 in ~\cite{ChronGear1989}), all from moments of the residual $r_i$.
In particular, for the matrix $W_i$, we have the following relation: $
W_i=M_i-c_i\beta_{i-1}$,
where $c_i = (c_i^j)_{j=1, \dots, s}$ and $\beta_{i-1} = (\beta_{i-1}^j)_{j=1, \dots, s}$ are $s \times s$ matrices whose columns are the vectors $c_i^j$ and $\beta_{i-1}^j$, respectively, and $M_i$ is the $2s \times 2s$ matrix of moments of $r_i$ with respect to $AB^{-1}$ and the $B$-inner product, which is defined as:
\begin{equation}
\label{mommat}
M_i=(\mu_i^{j+k}), \ \ 0 \leq j,k \leq s,
\end{equation}
where $\mu_i^{j+k}=((AB^{-1})^{j+k}r_i)^TB^{-1}r_i=r_i^T(B^{-1}A)^{j+k}B^{-1}r_i$.
The matrix $M_i$ is spd as long as the following $s$ vectors:
$\{ B^{-1}r_i, (B^{-1}A)B^{-1}r_i, \ldots, (B^{-1}A)^{s-1}B^{-1}r_i \}$
form a basis for the $s$-dimensional Krylov subspace $K_s(r_i)$, meaning that they are linearly independent. However, it is important to note that numerical instabilities, due to round-off errors from inexact arithmetic, may cause a loss of orthogonality in the monomial basis vectors as $s$ increases or as the matrix dimension $n$ becomes large, as pointed out in ~\cite{ChronGear1989}.
Finally, we summarize the preconditioned $s$-step CG algorithm, focusing on the main computational steps. Note that we use Matlab-like notation to indicate vector and matrix entry index intervals.
{\small 
\begin{algorithm}
  \caption{PCGs(A, b, B, x0, s, rtol, kmax) (Part 1)}
  \label{alg:pcgs}
  \begin{algorithmic}[1]
  \State $x=x_0 \ \ r=b-A x_0$
  \State $P1=0$
    \For{$k=0,1,...,kmax$}
       \If{($k\bmod 2 == 0$)}
       \State $P2=pMPK(A,r,s,B)$
       \State $vm = cvm(P2,r,s)$
        \If{($k == 0$)}
	\For{$i=1:s$}
            \For{$j=1:s$}
            \State $W(i,j)=vm(i+j)$ 
         \EndFor
        \EndFor
        \State $\beta=0$
	\State Solve for $\alpha$: $W \alpha = vm(1:s)^T$
    \Else
        \State Compute $\alpha,\beta, W$ using $SW(\alpha,\beta,W,vm,s)$
    \EndIf
        \State $P2(:,1:s)=P2(:,1:s)+P1(:,1:s)\beta$
        \State $x=x+P2(:,1:s) \alpha$
	   \State $P2(:,s+1:2s)=P2(:,s+1:2s)+P1(:,s+1:2s) \beta$
	   \State $r=r-P2(:,s+1:2s)\alpha$
    \EndIf
%
  \If{($k\bmod 2 == 1$)}
        \State $P1=pMPK(A,r,s,B)$
        \State $vm = cvm(P1,r,s)$
        \State Compute $\alpha,\beta,W$ using $SW(\alpha,\beta,W,vm,s)$
            \State $P1(:,1:s)=P1(:,1:s)+P2(:,1:s)\beta$
            \State $x=x+P1(:,1:s) \alpha$
	   \State $P1(:,s+1:2s)=P1(:,s+1:2s)+P2(:,s+1:2s)\beta$
	   \State $r=r-P1(:,s+1:2s)\alpha$
    \EndIf
    \If{($\|r\| < rtol$)}
        \State return $x$
    \EndIf
    \EndFor
    \State Return error, algorithm did not converge within $kmax$ iterations
  \end{algorithmic}
\end{algorithm}
}
{\small
\begin{algorithm}
  \caption{Preconditioned MP kernel $pMPK(A,x,s,B)$}
  \label{alg:pMPcg}
  \begin{algorithmic}[1]
  \State $P(:,1)=B^{-1}x$
  \State $pin=1, pout=s+1$
    \For{j=2:2s}
   	\If{($j\bmod 2 == 0$)}
   	\State $P(:,pout)=A P(:,pin)$ 
   	\Else
   	\State $P(:,pout)=B^{-1} P(:,pin)$ 
        \EndIf
        \State $ptemp=pin, pin=pout, pout=ptemp+1$
    \EndFor
    \State return $P$
  \end{algorithmic}
\end{algorithm}
\begin{algorithm}
  \caption{$cvm(P,r,s)$}
  \label{alg:cvm}
  \begin{algorithmic}[1]
        \State $vm(1:s)=P(:,1:s)^Tr$ 
        \State $vm(s+1:2s)=P(:,1:s)^TP(:,2s)$ 
        \State return $vm$
    \end{algorithmic}
\end{algorithm}
\begin{algorithm}
  \caption{$SW(\alpha,\beta, W,vm,s)$}
  \label{alg:pSWcg}
  \begin{algorithmic}[1]
	\State $rhs=zeros(2s-1,1)$
	\For{$i=1:s$}
	\State $rhs(s-1+i)=vm(i)$
    	\For{$j=1:i-1$}
    	\State $rhs(s-1+i)=rhs(s-1+i)+rhs(s-1+j)\alpha(s-i+j)$
    	\EndFor
    	\State $rhs(s-1+i)= -rhs(s-1+i)/\alpha(s)$
    	\EndFor
 	\For{$i=1:s$}
    	\For{$j=1:i-1$}
    	\State $b1(i,j)= -rhs(i+j-1)$
    	\EndFor
    	\EndFor
	\State Solve for $\beta$: $W\beta = b1$
 	\For{$i=1:s$}
            \For{$j=1:s$}
            \State $W(i,j)=vm(i+j)$ 
         \EndFor
        \EndFor
	\State $W=W-b1\beta$  	
	\State Solve for $\alpha$: $W\alpha = vm(1:s)^T$
    \State return $\alpha,\beta,W$
  \end{algorithmic}
\end{algorithm}
}

Algorithm \ref{alg:pcgs}, at each iteration, requires the application of the so-called preconditioned Matrix-Power (MP) kernel, as described in Algorithm \ref{alg:pMPcg}. In this kernel, a preconditioner application and an SpMV operation are alternately applied, with a total of $s$ preconditioner applications and $s$ SpMV operations. We highlight that in Algorithm \ref{alg:pMPcg}, custom memory access is employed for the alternately computed vectors $B^{-1}x$ and $Ax$. This memory access strategy enables the use of contiguous memory blocks for multiple vector updates in Algorithm \ref{alg:pcgs}.
Furthermore, Algorithm \ref{alg:cvm} simultaneously computes $2s$ dot products, which contributes to improving computational efficiency. Algorithm \ref{alg:pSWcg} performs $\sum_{i=1}^s i = \mathcal{O}(s^2 / 2)$ floating-point operations to update a vector and solve two linear systems of dimension $s$. These operations introduce a computational overhead, but they can typically be locally executed in parallel, and this overhead is generally negligible or amortized by the $\mathcal{O}(s)$ reduction in global synchronization steps when compared to the classic PCG formulation, as seen in Algorithm \ref{alg:pcghs1}.
Finally, $s$ axpy operations (which can be expressed in terms of $s$-dimensional matrix row updates) are computed.

\section{Multi-GPU Design and Implementation Issues}
\label{design}

Our design patterns for implementing Algorithm \ref{alg:pcgs} utilize GPU accelerators for all major computational components, optimizing data and thread locality. Additionally, they minimize data transfers between CPUs and GPUs within a node, as well as across multiple nodes. A key objective was to reuse data structures and functionalities from the existing software library {\em BootCMatchGX}, facilitating seamless integration of the new solver into the framework and leveraging prior developments. Specifically, we employed a contiguous row-block decomposition of the system matrix, stored locally using the Compressed Sparse Row (CSR) format, and implemented an efficient hybrid MPI-CUDA Sparse Matrix-Vector (SpMV) multiplication. This implementation overlaps computations on local matrix rows assigned to GPUs with asynchronous point-to-point communications of the matrix halo (or boundary) data. This operation is essential for executing the preconditioned Matrix Powers Kernel (MPK), as detailed in Algorithm \ref{alg:pMPcg}.
For more information on BootCMatchGX, you can visit its GitHub repository.

We have developed a specialized GPU kernel for Algorithm \ref{alg:cvm} to compute the $2s$ moments of the current residual efficiently. This kernel aggregates the dot products required for orthogonalizing Krylov basis vectors, necessitating only a single global reduction to derive the $2s$ scalars that define the matrix $W$ and the right-hand side vectors essential for obtaining the $s$-dimensional vector $\alpha$ and the $s$-dimensional matrix $\beta$. Our GPU kernel is designed to maximize data locality by reusing data already present in GPU registers, thereby avoiding unnecessary read and write operations to and from GPU memory. Since the global reduction involves MPI tasks running on CPU cores, all components of the $\alpha$ vector and entries of the $\beta$ matrix are accessible to all distributed MPI tasks. Consequently, the solution of the small linear system at line 14 of Algorithm \ref{alg:pcgs}, as well as computations in Algorithm \ref{alg:pSWcg}, are replicated by each MPI task and executed on CPU cores.
Matrix and vector updating recurrences at lines $18-19$ and $21-22$, as well as at lines $31-32$ and $34-35$, are purely local and can be efficiently performed on GPUs by leveraging their architectural features. Notably, this scenario allows immediate utilization of freshly computed data (e.g., results from matrix updates) for subsequent operations (e.g., vector updates) while the data resides in GPU registers. Instead of sequentially employing a combination of {\tt cublasDgemm} followed by {\tt cublasDgemv}—which would involve writing the matrix to GPU memory and reading it back—we implemented a custom {\tt Dgemmv} kernel that computes the final result directly. This approach enhances efficiency by reducing memory access overhead. In our custom kernel, the matrix product is optimized by loading all elements of the second matrix into the GPU's shared memory. This strategy is feasible for substantial values of $s$; for instance, in our experiments, $s$ can be as large as 90, fitting comfortably within the shared memory capacities of current GPUs. For smaller values of $s$ (e.g., $s \leq 5$), each matrix element of the result is computed by a single thread rather than a warp (a group of 32 threads), as this configuration has proven more efficient in our tests.
 While Algorithm \ref{alg:cvm} could be implemented using a sequence of {\tt cublasDgemv} primitives, our alternative approach—employing a custom multi-dot product within a single CUDA kernel—has demonstrated superior efficiency. Although these components may have a relatively modest impact on the overall execution time, our custom kernels have shown a performance improvement factor of 2 to 3 compared to combinations of {\tt cublas} primitives, depending on problem size and the value of $s$.

 \subsection{Parallel Preconditioner}
\label{prec}

Reducing communication and increasing computational granularity in a preconditioned Conjugate Gradient (CG) method necessitates a judicious selection of the preconditioner 
$B$. The goal is to balance convergence efficiency—ensuring the iteration count remains nearly constant regardless of problem size—with implementation scalability, achieved through a high degree of parallelism in applying the preconditioner.
A primary objective in developing {\em BootCMatchGX} was to design preconditioners whose application depends solely on the fundamental operation used in Krylov Subspace Methods (KSMs): the Sparse Matrix-Vector multiplication (SpMV). As previously discussed, SpMV involves only local communications, and {\em BootCMatchGX} provides a highly efficient hybrid MPI-CUDA implementation for this operation.
The preconditioner applied in Section~\ref{weak} is an Algebraic MultiGrid (AMG) method, specifically designed to be scalable in both convergence properties and parallelism. This design addresses the challenges of utilizing clusters of GPU-accelerated computing nodes for solving sparse linear systems at extreme scales. The operator $B^{-1}$ is recursively defined through several weighted Jacobi iterations on a hierarchy of coarser sparse matrices derived from the original system matrix. Consequently, the primary computational component is the Sparse Matrix-Vector (SpMV) product, applied to both the system matrix and a sequence of progressively smaller, though denser, matrices. Details on the AMG preconditioner and on the design patterns applied for its implementation on a multi-GPU clusters are in~\cite{parco2020,ieee2023}.

\section{Related Work}
\label{relwork}

The use of $s$-step variants of preconditioned Krylov Subspace Methods (KSMs) is highly appealing due to their potential to significantly reduce communication costs while increasing the computational granularity of linear algebra operations. These methods achieve this by performing vector orthogonalization in blocks of dimension $s$.
The foundational work by Chronopoulos and Gear in 1989~\cite{ChronGear1989p,ChronGear1989} introduced $s$-step methods for symmetric matrices, paving the way for extensions to unsymmetric and indefinite matrices. Carson et al.~\cite{carsonacta,carsonthesis} expanded on these methods, addressing numerical instabilities by proposing alternative polynomial bases for Krylov subspaces, such as Newton or Chebyshev bases, and adaptive strategies to select optimal $s$ values. Their innovations include:
i) residual replacement, i.e., recomputing residuals directly to improve stability and reduce errors due to recursive updates;
ii) mixed precision techniques, i.e., using higher precision for momentum-related computations, as demonstrated in multi-GPU implementations~\cite{carsonichitaro2022}.
Alternative synchronization reduction strategies, such as pipelining~\cite{GHYSELS2014}, enable overlapping SpMV computations with global communication steps. Pipelining has been integrated into MPI-based $s$-step CG solvers, as shown in~\cite{manasi2021}.
Some efforts have been devoted to the implementation of a Communication-Avoiding (CA-)MP kernel to compute the Krylov basis. This is fairly simple for some simple sparsity patterns of the system matrix, although it is not clear whether it is truly efficient for hybrid systems, in the version without preconditioning, but very challenging if a reliable preconditioner has to be applied (as the preconditioner itself requires communication). Therefore, it appears also clear why the design and the implementation of communication-reduced/avoiding preconditioners to be coupled to CA-KSM is a critical aspect~\cite{grigori2015,ichitaro2014,MGDD2021}.
Experimental studies highlight the practical impact of $s$-step CG methods with communication reduction techniques. On the K Japanese supercomputer, for instance, preconditioned $s$-step CG combined with a block-Jacobi preconditioner and CA-MP kernel demonstrated improved performance~\cite{Kcomputer2016}.

\section{Numerical Results}
\label{results}

In this section, we discuss some results of our multi-GPU implementation of Algorithm \ref{alg:pcgs} obtained on the Alex cluster operated by the Erlangen National High-Performance Computing Center (NHR) of the Friedrich-Alexander-University (FAU). That cluster is equipped with nodes having 2 AMD EPYC 7713 (Milan) processors (64 cores per chip) and 8 Nvidia A100 (each with 40 GB HBM2). Two HDR200 Infiniband HCAs connect the nodes. For our experiments, we used up to $8$ nodes for a total of $64$ GPUs. Nvidia Cuda (rel. 11.6) and OpenMPI are used on the top of the GNU toolchain.

In the following, we analyze both strong and weak scalability of the software module in solving linear systems arising from the discretization of the Poisson equation on a 3D cube, with Dirichlet boundary conditions and right-hand side equal to $1$.
We apply a classic five-point finite-difference discretization on a mesh with $250^3 \approx 15.6 \times 10^6$ points, corresponding to the full load for a single GPU. Solver procedures stopped when the Euclidean norm of the relative residual reached the tolerance $rtol = 10^{-6}$ or the number of iterations reached a predefined value $itmax = 2000$.

Our first aim is to compare results obtained by using the original {\em BootCMatch} solver (PFCG as in Algorithm \ref{alg:pfcgnn}), which is already a communication-reduced version of the original method, versus the new PCGs described in Algorithm \ref{alg:pcgs}. 

\subsection{Strong Scalability}
\label{strong}

In the following we discuss results obtained when no preconditioning is applied; this formally corresponds to use the identity matrix as preconditioner and it is interesting because the number of iterations, per variant of the method, is always the same regardless of the number of GPUs. The main aim is to analyze the impact of reducing global synchronization steps and increasing computational granularity when Algorithm \ref{alg:pcgs} is applied for increasing values of $s$. We use values of $s$ from $1$ to $5$, for which we did not encounter any numerical instability during the solution process. In this section we analyze the strong scalability behaviour of the solver, that is we fixed the size of the problem to $n=250^3$ and increased number of GPUs up to $64$. 
%

%

\begin{figure}[h!]
\includegraphics[width=0.5\textwidth]{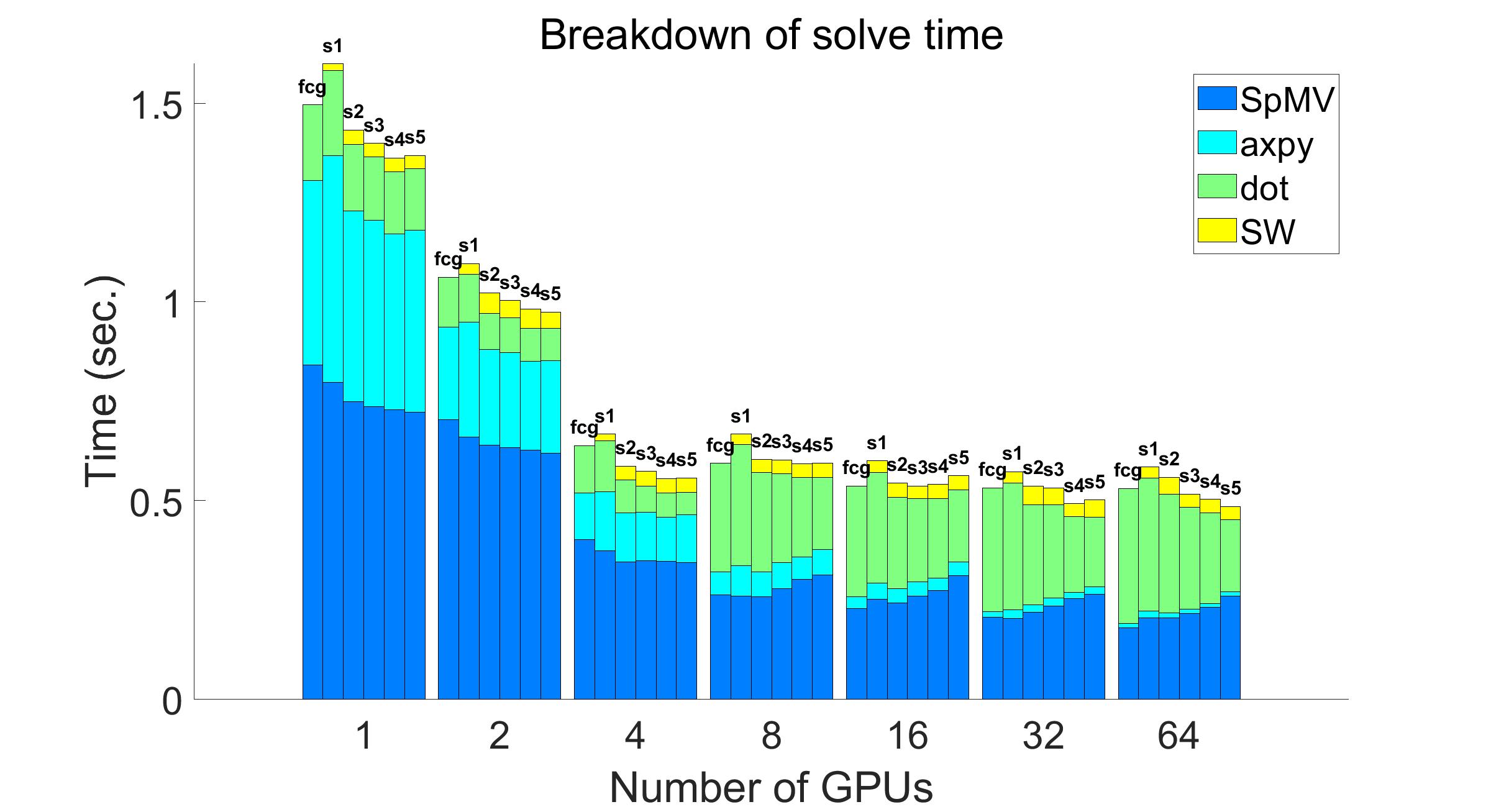}
\caption{Strong scalability: breakdown of solve time when no preconditioner is applied. \label{fig:breaknop}}
\end{figure}

In Figure \ref{fig:breaknop}, we present a breakdown of the solve time across the various computational building blocks, identified by different colors in the vertical bars, as detailed in the legend. Each group of bars includes the first bar corresponding to Algorithm \ref{alg:pfcgnn}, referred to as \verb|fcg| for brevity, and five additional bars corresponding to Algorithm \ref{alg:pcgs} for increasing values of 
$s$, labeled by the respective value of $s$. Note that the green color represents local dot computations and data communication for the global sum.

As expected, the most computationally expensive building block is the SpMV computation. Its cost decreases rapidly up to $8$ GPUs and continues to decline more slowly as the number of GPUs increases to $64$. This trend results from the decreasing computation-to-communication ratio, which occurs when the total problem size remains fixed while the number of GPUs increases. Note that the oscillations observed among the different solver versions—despite the same number of SpMVs being performed—remain a small fraction of the total execution time.

A key observation concerns the behavior of the axpy and dot computation blocks. As anticipated, axpy operations significantly impact the runtime on a single GPU, but their cost decreases rapidly as the number of GPUs increases. This is because axpy operations are embarrassingly parallel, and the fine-tuned kernel designed for efficient data reuse at the GPU level delivers substantial benefits. In contrast, dot computations exhibit the opposite behavior: their global communication costs increase as GPUs increase. In fact, for both PFCG and PCGs, when 
$s$ is small and the number of GPUs exceeds $8$, the costs of dot computations can become dominant, even surpassing those of SpMV operations.
It is worth noting that we observe a significant increase in the global reduction required for dot computations on a single node when going from $4$ to $8$ GPUs. This is confirmed by a separate benchmark of the \verb|MPI_Allreduce| primitive (that is used in our implementation of the dot), applied to a $64$-byte scalar variable, with MPI tasks ranging from $1$ to $32$.
Finally, we observe that the overhead computations (SW) in the application of the PCGs incur an almost constant cost, with only small increases as 
$s$ grows. However, this overhead is generally well amortized as $s$ increases, particularly on $64$ GPUs, due to the reduction in the cost of dot computations.

\begin{figure}[h!]
\includegraphics[width=0.5\textwidth]{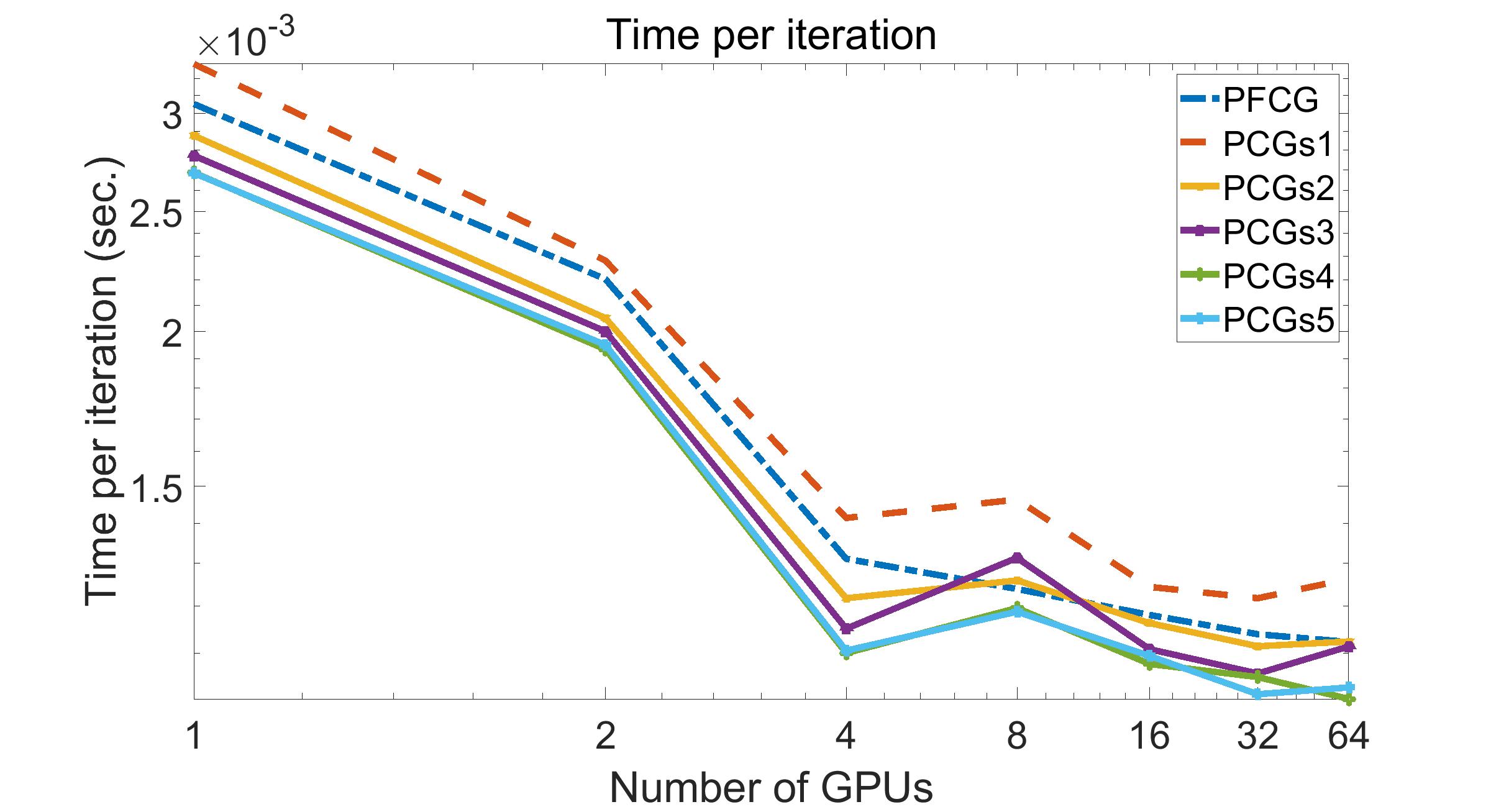}
\caption{Strong scalability: solve time per iteration when no preconditioner is applied. \label{fig:t4itnop}}
\end{figure}

To better understand the behavior of the various solver versions, Figure \ref{fig:t4itnop} presents the execution times per iteration of the iterative process. For a fair comparison, the cost of an iteration of Algorithm \ref{alg:pcgs} is scaled by $s$, since one iteration of PCGs corresponds to $s$ iterations of PFCG in terms of solution improvement. It is clear that, for values of $s$ greater than $2$, PCGs exhibit a smaller time per iteration, which generally decreases as the number of GPUs increases. However, on $8$ GPUs, a slight increase is observed, as the reduction in SpMV computation is insufficient to offset the significant increase in dot computation. This behavior highlights the importance of minimizing global synchronizations and communication in KSMs.

\subsection{Weak Scalability}
\label{weak}
In this section, we discuss the performance results of the solvers in a weak scaling setting, where the problem size increases linearly with the number of GPUs. As in the previous section, we fix the problem size on a single GPU at $n=250^3$, so that we solve a system with 1 billion ($n=10^9$) unknowns on 64 GPUs.

\begin{figure}[h!]
\includegraphics[width=0.5\textwidth]{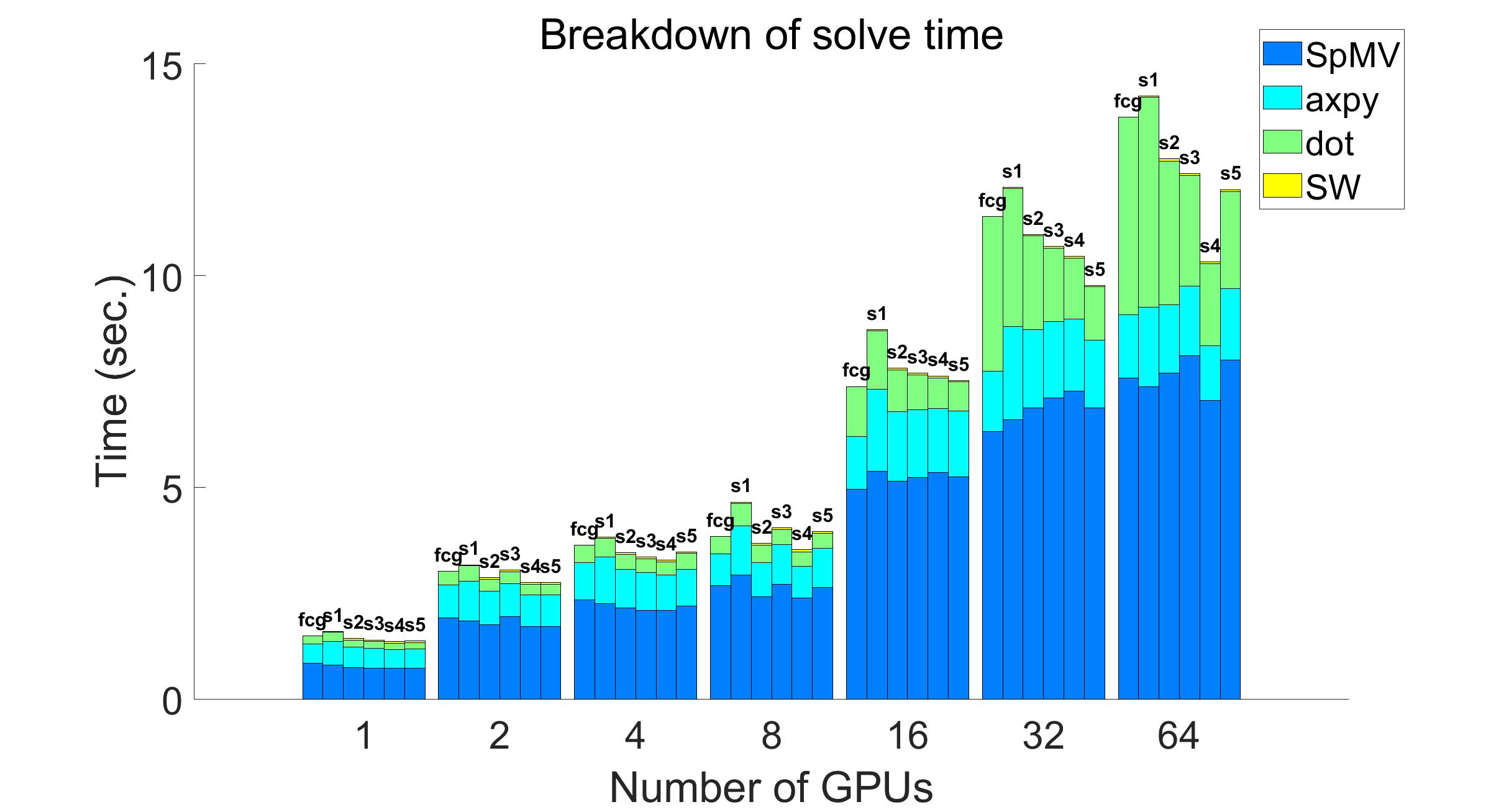}
\caption{Weak Scalability: breakdown of solve time when no preconditioner is applied. \label{fig:breaknopws}}
\end{figure}

The breakdown of the solve time in Figure \ref{fig:breaknopws} shows that, as expected for this type of computation, the total solve time increases with the problem size. This is partly because the number of iterations required by the various solver versions, when no preconditioner is applied, increases significantly to achieve the desired convergence. This increase is caused by the growing maximum-to-minimum eigenvalue ratio as the problem dimension expands.

Specifically, the number of iterations for PFCG and PCGs with $s=1$ ranges from $514$ on $1$ GPU to $2101$ on $64$ GPUs. Interestingly, while axpy operations and SW computations show minimal increases in time as the number of GPUs increases, the largest cost is due to the SpMV computation. The ratio of this cost increase, from $1$ GPU to $64$ GPUs, is approximately $10$. In contrast, the cost of dot computations rises significantly with the number of GPUs. For PFCG and PCGs with 
$s=1$, the increase ratio is around $20$, while for PCGs with 
$s=5$, the ratio is reduced to about $11$, demonstrating the benefits of reducing communication and synchronization steps.

\begin{figure}[h!]
\includegraphics[width=0.5\textwidth]{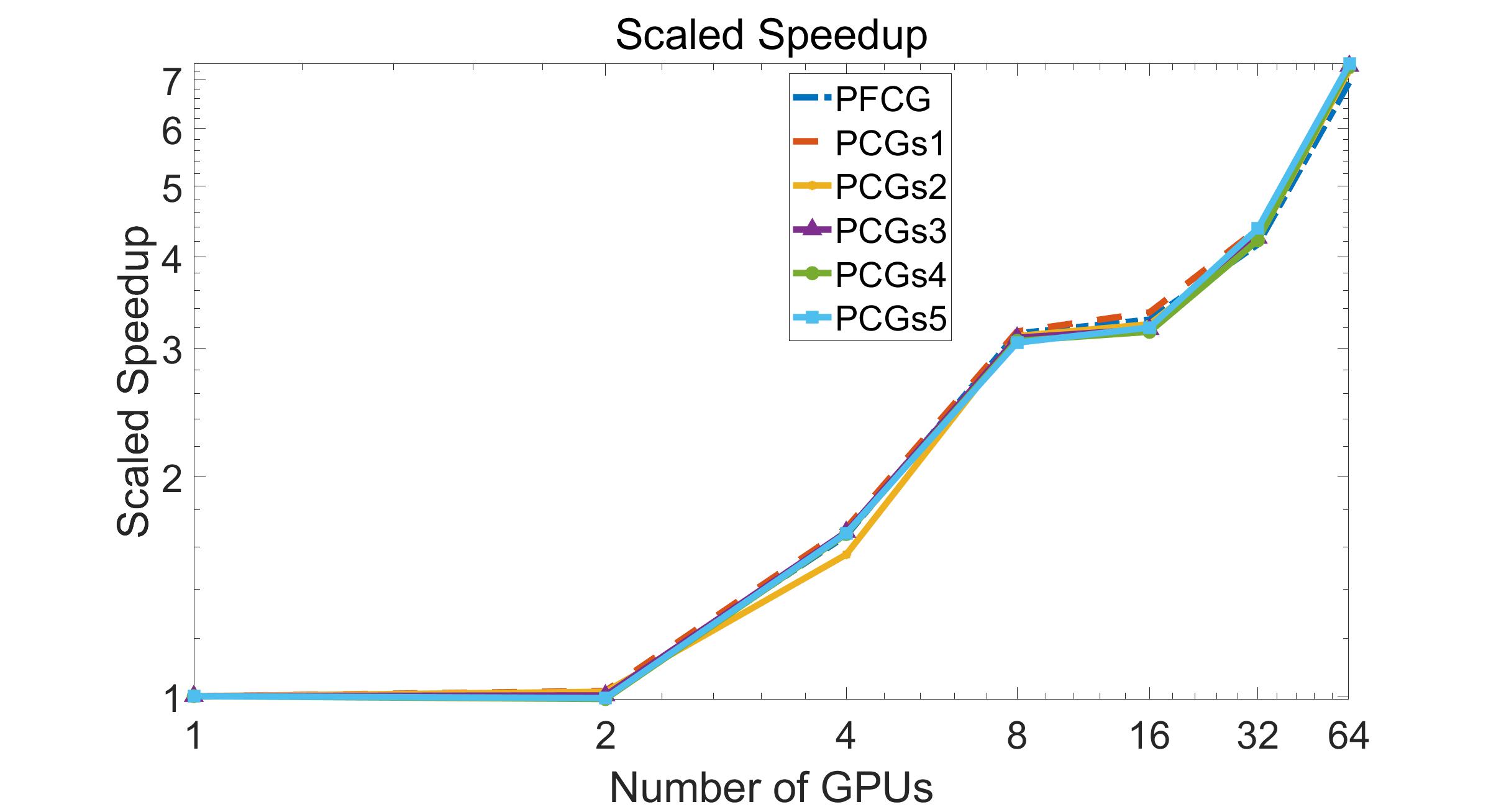}
\caption{Scaled Speedup of the solve time when no preconditioner is applied. \label{fig:speedupnopws}
}
\end{figure}

The best solve time with an increasing number of GPUs is achieved by PCGs, especially for larger values of $s$. However, it is important to note that for $s=4$ and $64$ GPUs, numerical instabilities occur, causing the solver to exit before the desired accuracy is reached.
For completeness, Figure \ref{fig:speedupnopws} presents the scaled speedup, defined as $SP=\frac{T_1(n) \times p}{T_p(n \times p)}$, where $T_1(n)$ is the solve time for a problem size $n$ on $1$ GPU and $T_p(n \times p)$ is the solve time for a problem of size $n \times p$ on $p$ GPUs. We observe that the scaled speedup values are similar for all variants of the solvers. However, as the number of GPUs increases, the best values are generally achieved for the largest values of 
$s$. For instance, on $64$ GPUs, the speedup ranges from about $7$ for PFCG to approximately $7.36$ for PCGs. In light of the behavior of the ratio factor in dot computations, we expect this gap to widen as the number of GPUs increases.

%
%
%
\begin{figure}[h!]
\includegraphics[width=0.5\textwidth]{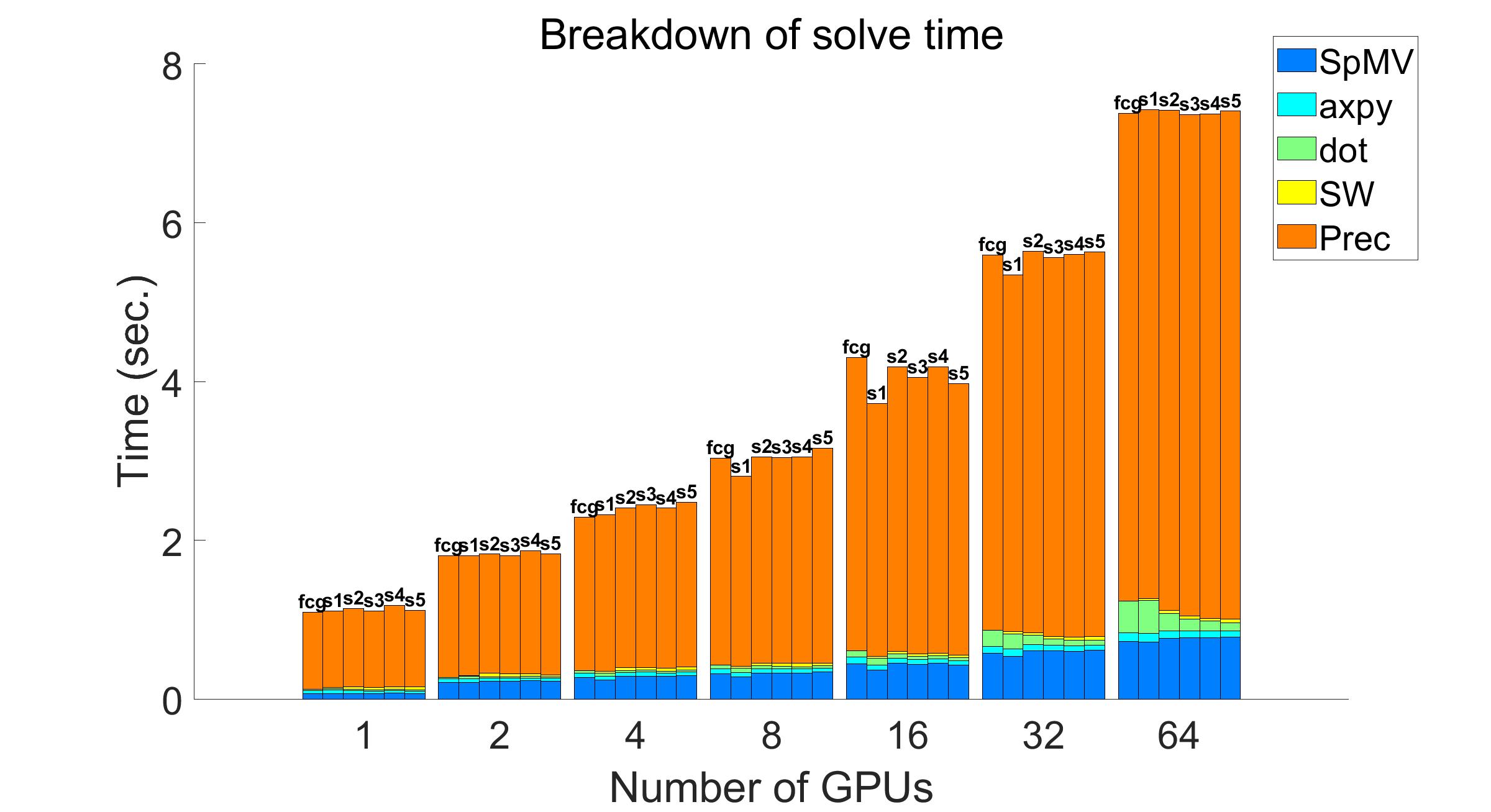}
\caption{Weak Scalability: breakdown of solve time when AMG preconditioner is applied. \label{fig:breakamgpws}}
\end{figure}
%
%
\begin{figure}[h!]
\includegraphics[width=0.5\textwidth]{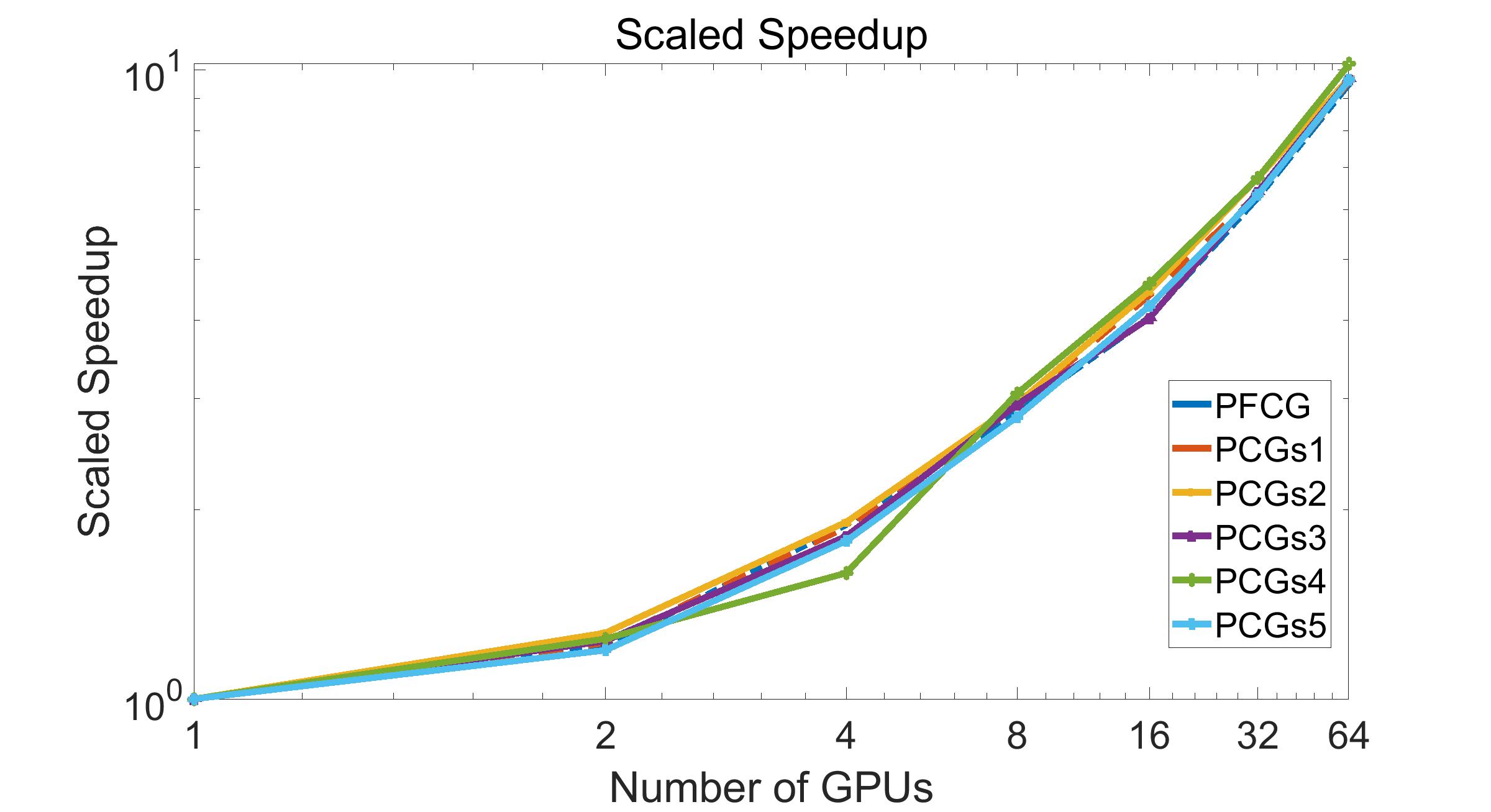}
\caption{Scaled Speedup of the solve time when AMG preconditioner is applied. \label{fig:speedupamgpws}}
\end{figure}

In Figure \ref{fig:breakamgpws}, we show the breakdown of the solve time when the AMG preconditioner is applied. In this case, no instability is observed for any of the $s$-values used in the PCGs variant.

As expected, the application of the preconditioner is the most computationally expensive operation. However, it leads to the best solve times as the number of GPUs increases, primarily due to the slower increase in the number of iterations, which range from $45$ to $119$. On the other hand, the reduction in the cost of dot product computations, which decreases from approximately $21$ to $5$ as 
$s$ increases for PCGs, is nearly offset by the reduction in the number of iterations. Nevertheless, this reduction is expected to have a larger impact as the number of GPUs and the problem size increase, particularly for larger values of $s$. Finally, as shown in Figure \ref{fig:speedupamgpws}, we observe that in this case, the scaled speedup is approximately $10$ on $64$ GPUs for all variants of the solver.
\section{Concluding Remarks}
\label{conclusions}

In this paper, we presented the design principles, implementation strategies, and performance results of communication-reduced variants of the widely used preconditioned Conjugate Gradient (CG) method for solving large and sparse linear systems on clusters of GPU-accelerated computing nodes.

The primary innovations lie in the design of software modules that leverage optimized data communication and access patterns, as well as data reuse strategies, to fully exploit the high throughput of GPUs in large, heterogeneous clusters. These advancements are made available through a freely accessible scientific library that can be seamlessly integrated with application codes.

Future work will focus on re-designing algorithmic and implementation approaches to enhance numerical stability. Additionally, the development of further communication-reduced parallel preconditioners—capable of balancing convergence efficiency with high degrees of parallelism—remains a key area of investigation.

\section{Acknowledgements}
\label{acks}

The authors gratefully acknowledge the scientific support and HPC resources provided by the Erlangen National High Performance Computing Center (NHR@FAU) of the Friedrich-Alexander-Universit\"{a}t Erlangen-N\"{u}rnberg (FAU). The hardware is funded by the German Research Foundation (DFG) and  federal and Bavarian state authorities through the NHR program.

\bibliographystyle{IEEEtran}
\bibliography{bibliopdp2025}

\end{document}